\documentclass[11pt]{article}
\usepackage{amssymb,amsfonts,amsmath,amsthm}
\usepackage{epsfig}
\newtheorem{thm}{Theorem}[section]

\newtheorem{ex}[thm]{Example}
\newtheorem{alg}[thm]{Algorithm}

\makeatletter \@addtoreset{equation}{section}

\begin{document}

\begin{center}
{\Large\bf Converging to Gosper's Algorithm}

\vskip 5mm

William Y. C. Chen$^1$, Peter Paule$^2$, and Husam L. Saad$^3$

 $^{1,3}$Center for Combinatorics, LPMC-TJKLC\\
 Nankai University, Tianjin 300071, P. R. China

 \vskip 1mm

 $^{2}$Research Institute for Symbolic Computation,\\
J. Kepler University, A-4040 Linz, Austria

$^1$chen@nankai.edu.cn, \ $^2$Peter.Paule@risc.uni-linz.ac.at,
 \
$^3$hus6274@hotmail.com

\end{center}

\vskip 3mm

\centerline{\bf Abstract}

 Given two polynomials, we find a
convergence property of the GCD of the rising factorial and the
falling factorial. Based on this property, we present a unified
approach to computing the universal denominators  as  given by
Gosper's algorithm and Abramov's algorithm for finding rational
solutions to linear difference equations with polynomial
coefficients.

\vskip 3mm



\vskip 3mm

\section{Introduction}\label{sec1}

Let $\mathbb N$ be the set of nonnegative integers, $\mathbb K$ be
a field of characteristic zero, $\mathbb K(n)$ be the field of
rational functions over $\mathbb K$, and $\mathbb K[n]$ be the
ring of polynomials over $\mathbb K$. We assume that subject to
normalization the gcd (greatest common divisor) of two polynomials
always takes a value as a monic polynomial, namely, polynomials
with the leading coefficient being 1. Recall that a nonzero term
${t_{n}}$ is called a hypergeometric term over $\mathbb K$ if
there exists a rational function $r\in\mathbb K(n)$ such that
\begin{equation*}
\frac{t_{n+1}}{t_{n}}=r(n).
\end{equation*}
\noindent If $r(n)=a(n)/b(n)$, where $a(n), b(n)\in\mathbb K[n]$, then the
function $a(n)/b(n)$ is called a rational representation of the rational
function $r(n)$. If $\mathrm{gcd}(a(n),b(n))=1$ holds, then $a(n)/b(n)$ is
called a reduced rational representation of $r(n)$.

Gosper's algorithm \cite{Gos78} (also see \cite{GKP89, Koe98,
PWZ96,Zei90a, Zei90b, Zei91,WZ92}) has been extensively studied
and widely used to prove hypergeometric identities. Given a
hypergeometric term $t_{n}$, Gosper's algorithm is a procedure to
find a hypergeometric term $z_{n}$ satisfying
\begin{equation}\label{eq11}
z_{n+1}-z_{n}=t_{n},
\end{equation}
if it exists, or confirm the nonexistence of any solution of
(\ref{eq11}). The key idea of Gosper's algorithm lies in a
representation of rational functions called Gosper representation;
i.e., writing the rational function $r(n)$ in the following form:
\begin{equation*}
r(n)=\frac{a(n)}{b(n)}\frac{c(n+1)}{c(n)},
\end{equation*}
where $a$, $b$ and $c$ are polynomials over $\mathbb K$ and
\begin{equation*}
\mathrm{gcd}(a(n),b(n+h))=1 \;\;\mbox{for all }  h\in\mathbb N.
\end{equation*}
Petkov\v{s}ek \cite{Pet92} has realized that a Gosper representation
becomes unique, which is called the Gosper-Petkov\v{s}ek
representation, or GP representation, for short, if we further
require that $b$, $c$ are monic polynomials such that
\begin{equation*}
\mathrm{gcd}(a(n),c(n))=\mathrm{gcd}(b(n),c(n+1))=1.
\end{equation*}

In the same paper, Petkov\v{s}ek also gave an algorithm to compute
GP representations; subsequently we will call it  the ``GP
algorithm''. In \cite{PS95}, Paule and Strehl gave a derivation of
Gosper's algorithm by using the GP representation. In \cite{Pau95},
equipped with the Greatest Factorial Factorization (GFF), Paule
presented a new approach to indefinite hypergeometric summation
which leads to the same algorithm as Gosper's, but in a new setting.
In \cite{LPS93}, Lison\v{e}k and {\it et al.,} gave a detailed study
of the degree setting for Gosper's algorithm.

Finding rational solutions is important in computer algebra
because many problems can be reduced to rational solutions. For
example, we may consider the generalization of Gosper's algorithm.
Given a linear difference equation
\begin{equation}\label{eq12}
\sum_{m=0}^{d}p_{m}(n)y(n+m)=p(n),
\end{equation}
where $p_{0}(n), p_{1}(n), \ldots, p_{d}(n), p(n)\in\mathbb K[n]$
are given polynomials such that $p_{0}(n)\neq0$, $p_{d}(n)\neq0$, a
polynomial $g(n)\in\mathbb K[n]$ is called a universal denominator
for (\ref{eq12}) if and only if for every solution $y(n)\in\mathbb
K(n)$ to (\ref{eq12}) there exists a $f(n)\in\mathbb K[n]$ such that
$y(n)=f(n)/g(n)$. Once a universal denominator is found, then it is
easy to find the rational solutions of the linear difference
equation (\ref{eq12}) by finding the polynomial solutions using the
techniques in \cite{Abr89, ABP95, Pet92}. Abramov \cite{Abr89}
developed an algorithm to find a universal denominator of
(\ref{eq12}) which relies on all the coefficients $p_{0}(n),
p_{1}(n), \ldots, p_{d}(n), p(n)$. In \cite{Abr95}, an improved
version is given which requires only two coefficients $p_{0}(n)$ and
$p_{d}(n)$. Compared with the simplicity of the output of Abramov's
algorithm, the justification is quite involved. Recall that the
dispersion $\mathrm{dis}(a(n),b(n))$ of the polynomials
$a(n),b(n)\in\mathbb K[n]$ is the greatest nonnegative integer $k$
(if it exists) such that $a(n)$ and $b(n+k)$ have a nontrivial
common divisor, i.e.,
$$\mathrm{dis}(a,b)=\mathrm{max} \{k\in\mathbb N\ |\ \mathrm{deg\
gcd}(a(n), b(n+k))\geq1\}.$$ If $k$ does not exist then we set
$\mathrm{dis}(a,b)=-1$. Observe that $\mathrm{dis}(a(n),b(n))$ can
be computed as the largest nonnegative integer root of the
polynomial $R(h)\in\mathbb K[h]$ where $R(h)=Res_{n}(a(n),b(n+h))$.

The main result of this paper is the discovery of a convergence
property of the GCD of rising factorial of a polynomial $b(n)$ and
the falling factorial of another polynomial $a(n)$. By using the
limit of the GCD sequence, we may transform a rational difference
equation into a polynomial difference equation.  The convergence
argument yields a new and streamlined approach to the explicit
formula for Abramov's universal denominator. Note that this explicit
formula can be used to compute rational solutions of a linear
difference equation (\ref{eq12}). In addition, we derive Abramov's
universal denominator from Barkatou's explicit formula. The relation
between  Barkatou's approach and Abramov's algorithm has been
discussed in detail by Weixlbaumer \cite{Wei01}.

\section{The Convergence Property}\label{sec2}

The main idea of this paper is the following convergence property of
a GCD sequence. It turns out that this simple observation plays a
fundamental role in finding rational solutions of linear difference
equations, and it can be viewed as a unified approach to several
well-known algorithms.

\begin{thm}\label{thm21}
Let $a(n)$ and $b(n)$ be two nonzero polynomials in $n$ and let
\begin{equation}\label{eq21}
k_{0}=\mathrm{dis}(a(n-1),b(n))=\mathrm{max} \{k\in\mathbb N\ |\
\mathrm{deg\ gcd}(a(n-1), b(n+k))\geq1\}.
\end{equation}
Define
\begin{equation}\label{eq22}
G_{k}(n)=\mathrm{gcd}(b(n)b(n+1)\ldots b(n+k-1),a(n-1)a(n-2)\ldots
a(n-k)).
\end{equation}
Then the sequence $G_{1}(n), G_{2}(n), \ldots $ converges to
$G_{k_{0}+1}(n)$.
\end{thm}

\begin{proof} For all $k>k_0$ we have
\begin{equation}\label{eq23}
\mathrm{gcd}(a(n-1), b(n+k)) =\cdots = \mathrm{gcd}(a(n-k-1),
b(n+k)) =1.
\end{equation}
Note that \[ G_{k+1}(n) = \mathrm{gcd}(b(n)b(n+1)\ldots
b(n+k),a(n-1)a(n-2)\ldots a(n-k-1) ) .\]  This implies that
\[ G_k(n) =
G_{k+1}(n), \] for all $k>k_0$. Moreover, one sees that once
(\ref{eq23}) is satisfied for $k>k_0$, it is also satisfied for
$k+1$. It follows that
\[ G_k(n)= G_{k+1}(n) = G_{k+2}(n) = \cdots ,\]
for all $k>k_{0}$, and this completes the proof.
\end{proof}

 Using the above convergence property of the sequence
$G_k(n)$, we are led to a simple approach to Gosper's algorithm
without resorting to the Gosper representation or GP representation
of rational functions. Given a hypergeometric term $t_{n}$ and
suppose that there exists a hypergeometric term $z_{n}$ satisfying
equation (\ref{eq11}), then by using (\ref{eq11}) we find
\begin{equation}\label{eq24}
r(n)y(n+1)-y(n)=1,
\end{equation}
where $r(n)=t_{n+1}/t_{n}$ and $y(n)=z_{n}/t_{n}$ are rational
functions of $n$, see \cite{PWZ96}.\\

\begin{thm}\label{thm22} Let $r(n)$ and $y(n)$ in equation
$(\ref{eq24})$ be in terms of their reduced rational
representations$:$
\begin{equation}\label{eq25}
r(n)=\frac{a(n)}{b(n)},\ y(n)=\frac{f(n)}{g(n)}.
\end{equation}
Then $$g(n)\ |\ G_{k_{0}+1}(n),$$ where $k_{0}$ and $G_{k}(n)$ are
defined in Theorem $\ref{thm21}$.
\end{thm}

\noindent {\it Proof.} Using (\ref{eq25}) in (\ref{eq24}) gives
\begin{equation}\label{eq26}
a(n)g(n)f(n+1)-b(n)g(n+1)f(n)=b(n)g(n)g(n+1).
\end{equation}
From the above relation, we immediately get that
$$g(n)\,|\,  b(n)g(n+1)\ \ \ \mathrm{and\ that}\ \ \ \ g(n+1) \ |
\ a(n)g(n).$$ Using these two relations repeatedly we obtain
\begin{eqnarray*}
&g(n)\mid b(n)b(n+1)\ldots b(n+k-1)g(n+k),\\[6pt]
&g(n)\mid a(n-1)a(n-2)\ldots a(n-k)g(n-k),
\end{eqnarray*}
for all $k\in\mathbb N$. Since $\mathbb K$ has characteristic
zero,
$$\mathrm{gcd}(g(n),g(n+k))=\mathrm{gcd}(g(n),g(n-k))=1,$$ for all
large enough $k$. It follows that
\begin{eqnarray}
& g(n) \mid b(n)b(n+1)\ldots b(n+k-1),  \label{eq27} \\[6pt]
& g(n) \mid a(n-1)a(n-2)\ldots a(n-k), \label{eq28}
\end{eqnarray}
for all large enough $k$. Therefore
\begin{equation*}
g(n)\mid G_{k}(n),
\end{equation*}
for all large enough $k$. The rest of the proof follows when $k$
goes to infinity in this equation and by Theorem \ref{thm21}. \qed

The next step is simply to set
\begin{equation}\label{eq29}
g(n)=G_{k_{o}+1}(n)
\end{equation}
in equation (\ref{eq26}) as in the GFF algorithm of Paule. If
equation $\mathrm{(\ref{eq26})}$ can be solved for $f\in\mathbb
K[n]$, then
\begin{equation*}
z_{n}=\frac{f(n)}{g(n)}\ t_{n}
\end{equation*}
is a hypergeometric solution of $\mathrm{(\ref{eq11})}$; Otherwise
no hypergeometric solution of $\mathrm{(\ref{eq11})}$ exists. Note
that the solution $f(n)$ may not be coprime to $g(n)$. However, it
is clear that this does not affect the solution of $y(n)$. Indeed,
the polynomials $f(n)$ and $g(n)$ can be recovered from the solution
of $y(n)$ after dividing the greatest common factors.

\begin{alg}\label{alg21}\hskip 0.1 cm\\
\noindent INPUT: $r(n)\in\mathbb K(n)$
such that $t_{n+1}/t_{n}=r(n)$ for large enough $n$ in $\mathbb N$.\\
OUTPUT: a hypergeometric solution $z_{n}$ of $\mathrm{(\ref{eq11})}$
if it exists, otherwise ``no hypergeometric solution of \eqref{eq11}
exists''.
\begin{itemize}
\item[$\mathrm {(1)}$] Decompose $r(n)$ into $a/b$ where $a$, $b$
are two relatively prime polynomials. \item[$\mathrm {(2)}$]
Compute $k_{0}$ as in $(\ref{eq21})$. \item[$\mathrm {(3)}$] If
$k_{0}\geq0$ then compute $g(n)=G_{k_{o}+1}(n)$, where $G_{k}(n)$
is defined as in $(\ref{eq22})$, otherwise $g(n)=1$.
\item[$\mathrm {(4)}$] If equation $\mathrm{(\ref{eq26})}$ can be
solved for $f\in\mathbb K[n]$ then return $z_{n}=\frac{f(n)}{g(n)}\
t_{n}$; Otherwise return ``no hypergeometric solution of
$\mathrm{(\ref{eq11})}$ exists''.
\end{itemize}
\end{alg}

Let us take an example from \cite{PWZ96}:

\begin{ex}\label{ex21} Let $t_{n}=(4n+1)\cdot\frac{n!}{(2n+1)!}$, then
$$r(n)=\frac{t_{n+1}}{t_{n}}=\frac{4n+5}{2(4n+1)(2n+3)}.$$ Hence $a(n)=4n+5$,
$b(n)=2(4n+1)(2n+3)$ and then $k_{0}=0$. Note that for all
$k>k_{0}$, equation $(\ref{eq23})$ is satisfied. From
$(\ref{eq29})$, $g(n)=n+\frac{1}{4}$. By $\mathrm{(\ref{eq26})}$,
$f(n)$ is a polynomial which satisfies
$$2f(n+1)-4(2n+3)f(n)=(2n+3)(4n+1).$$ The polynomial
$f(n)=-\frac{1}{2}(2n+1)$ is a solution of this equation. Therefore,
\[ z_{n}=\frac{f(n)}{g(n)}\ t_{n} =-2\ \frac{n!}{(2n)!}.\]
\end{ex}

We remark that the argument for the relations (\ref{eq27}) and
(\ref{eq28}) is used by Petkov\v{s}ek \cite{Pet92}. Moreover, the
products on the right hand sides of (\ref{eq27}) and (\ref{eq28})
can be written in the notation of rising or falling factorials as
introduced by Paule \cite{Pau95}.

\section{ Connections to Gosper's and Abramov's Algorithms}\label{sec3}

We will show how Theorem \ref{thm22} is related to Gosper's
algorithm and Abramov's algorithm for finding rational solutions of
linear difference equations with polynomial coefficients
\cite{Abr95}.

\noindent {\bf Abramov's Algorithm (general order $d$)}: Consider
the difference equation
\begin{equation}\label{eq31}
p_{d}(n)y(n+d)+\ldots +p_{0}(n)y(n)=p(n)
\end{equation}
with given $p_{0}(n), p_{1}(n), \ldots, p_{d}(n), p(n)\in\mathbb
K[n]$ such that $p_{0}$ and $p_{d}$ are nonzero. Abramov gave the
following algorithm to compute a universal denominator $G(n)$ for
(\ref{eq31}): Define
\begin{equation}\label{eq32}
N=\mathrm{dis}(p_{d}(n-d),p_{0}(n))=\mathrm{max}\{k\in\mathbb N\
|\ \mathrm{deg\ gcd}(p_{d}(n-d),p_{0}(n+k))\geq1\}.
\end{equation}
If $N=-1$ set $G(n)=1$, i.e., in this case all rational solutions
are polynomials. If $N\geq 0$, define
\begin{equation}\label{eq33}
A_{N+1}(n)=p_{d}(n-d),\ \ B_{N+1}(n)=p_{0}(n),
\end{equation}
and for $i=N$ down to $i=0$ do:
\begin{eqnarray}
d_{i}(n)&=&
\mathrm{gcd}(A_{i+1}(n),B_{i+1}(n+i)),\label{eq34}\\[6pt]
A_{i}(n)&=&\frac{A_{i+1}(n)}{d_{i}(n)}\ \ \ \mathrm{and}\ \ \
B_{i}(n)=\frac{B_{i+1}(n)}{d_{i}(n-i)}.\label{eq35}
\end{eqnarray}
If we use the notation of the falling factorial of a polynomial
introduced in \cite{Pau95} by
\[ [f(x)]^{\underline{k}} = f(x) f(x-1) \cdots f(x-k+1),\]
then Abramov's universal denominator of (\ref{eq31}) can be
written as
\begin{equation}\label{eq36}
G(n)=[d_{0}(n)]^{\underline{1}}[d_{1}(n)]^{\underline{2}}
\cdots[d_{N}(n)]^{\underline{N+1}}.
\end{equation}
There is also an explicit formula for Abramov's universal
denominator (\ref{eq36}), namely,
\begin{equation}\label{eq37}
G(n)=\mathrm{gcd}([p_{0}(n+N)]^{\underline{N+1}},[p_{d}(n-d)]^
{\underline{N+1}}).
\end{equation}
It can be seen that (\ref{eq36}) is equivalent to Theorem 3 in
Abramov-Petkov\v{s}ek-Ryabenko \cite{APR}. A generalized form of
(\ref{eq37}) can be found in Barkatou \cite{Bar99}. For
completeness, we give a proof of the fact that the presentations
(\ref{eq36}) and (\ref{eq37}) indeed coincide with each other. To
this end we will follow the survey of Weixlbaumer \cite{Wei01}.

 First of all, based on the definition of $N$ and the fact that
$A_{i}\, |\, A_{i+1}$ and $B_{i}\, |\, B_{i+1}$,  it follows  that
for $0\leq i\leq N+1$, we have
\begin{equation}\label{eq38}
\mathrm{gcd}(A_{i}(n),B_{i}(n+k))=1\ \ \hbox{for all}\ \ k>N.
\end{equation}
Moreover, from (\ref{eq34}) and (\ref{eq35}), for $0\leq i\leq N$
it follows that
\begin{equation}\label{eq39}
\mathrm{gcd}(A_{j}(n),B_{j}(n+i))=1\ \ \hbox{for}\ \  0\leq j\leq
i.
\end{equation}
Therefore,
\begin{eqnarray*}
G(n)&=& \mathrm{gcd}([p_{0}(n+N)]^{\underline{N+1}},[p_{d}(n-d)]^
{\underline{N+1}})\\[6pt]
&=& \mathrm{gcd}([B_{N+1}(n+N)]^{\underline{N+1}},[A_{N+1}(n)]^
{\underline{N+1}})\\[6pt]
&=&
\mathrm{gcd}\left(\left[\frac{B_{N+1}(n+N)}{d_{N}(n)}\right]^{\underline{N+1}},
\left[\frac{A_{N+1}(n)}{d_{N}(n)}\right]^{\underline{N+1}}\right)\cdot
[d_{N}(n)]^
{\underline{N+1}}\\[6pt]
&=& [d_{N}(n)]^ {\underline{N+1}}\cdot
\mathrm{gcd}([B_{N}(n+N)]^{\underline{N+1}},[A_{N}(n)]^
{\underline{N+1}}).
\end{eqnarray*}
Observe that
\begin{eqnarray*}
\lefteqn{\mathrm{gcd}(B_{N}(n+N),[A_{N}(n)]^ {\underline{N+1}})}
 \qquad \\[6pt]
&=&\mathrm{gcd}(B_{N}(n+N),A_{N}(n))\hskip 1.5 cm
\hbox{[by (\ref{eq38})]}\\[6pt]
&=&1,\hskip 5.3 cm \hbox{[by (\ref{eq39})]}
\end{eqnarray*}
and similarly,
\begin{equation*}
\mathrm{gcd}([B_{N}(n+N-1)]^{\underline{N}},A_{N}(n-N))=1.
\end{equation*}
Consequently,
\begin{equation*}
\mathrm{gcd}([B_{N}(n+N)]^{\underline{N+1}},[A_{N}(n)]^
{\underline{N+1}})=\mathrm{gcd}([B_{N}(n+N-1)]^{\underline{N}},[A_{N}(n)]^
{\underline{N}}).
\end{equation*}
In the same manner one can successively split off the factors
\[ [d_{N-1}(n)]^ {\underline{N}},\ldots,
[d_{0}(n)]^{\underline{1}}\] until one arrives at (\ref{eq36}),
which completes the proof of the equality of (\ref{eq36}) and
(\ref{eq37}).

Next we remark that Theorem \ref{thm22} follows from Abramov's
algorithm, strictly speaking, the universal denominator given by
Abbramov's algorithm.  Equation (\ref{eq24}) is equivalent to
\begin{equation}\label{eq310}
a(n)y(n+1)-b(n)y(n)=b(n),
\end{equation}
which is (\ref{eq31}) with $d=1$, $p_{1}(n)=a(n)$,
$p_{0}(n)=-b(n)$, and $p(n)=b(n)$. From (\ref{eq37}), Abramov's
algorithm gives the following  universal denominator of
(\ref{eq310})
\begin{equation}\label{eq311}
G(n)=\mathrm{gcd}([a(n-1)]^{\underline{N+1}},[b(n+N)]^
{\underline{N+1}}).
\end{equation}
where $N=k_{0}$ by (\ref{eq21}). Using (\ref{eq22}) we have that
$G(n)=G_{k_{0}+1}(n)$, hence Theorem \ref{thm22} determines the same
universal denominator as Abramov's algorithm.

Next we show that Abramov's algorithm delivers a Gosper
representation for $r(n)=t_{n+1}/t_{n}$ if $r(n)=a(n)/b(n)$ is the
reduced rational representation of $r(n)$. From (\ref{eq35}) we
obtain that $a(n-1)=A_{N+1}(n)=d_{N}(n)A_{N}(n)$, and by iteration,
\begin{equation}\label{eq312}
a(n)= d_{0}(n+1)\ d_{1}(n+1) \ldots d_{N}(n+1)\ A_{0}(n+1).
\end{equation}
Analogously, (\ref{eq35}) implies that
\begin{equation}\label{eq313}
b(n)=-d_{0}(n)\ d_{1}(n-1) \ldots d_{N}(n-N)\ B_{0}(n).
\end{equation}
Consequently, in view of representation (\ref{eq36}) for $G(n)$,
one obtains that
\begin{equation}\label{eq314}
\frac{a(n)}{b(n)}=-\frac{G(n+1)}{G(n)}\
\frac{A_{0}(n+1)}{B_{0}(n)}.
\end{equation}
Note that in (\ref{eq310}) w.l.o.g. we can assume that
$\mathrm{gcd}(a(n),b(n))=1$, which then implies
\begin{equation}\label{eq315}
\mathrm{gcd}(A_{0}(n+1),B_{0}(n))=1.
\end{equation}
But more is true. Namely, by (\ref{eq39})
\begin{equation*}
\mathrm{gcd}(A_{0}(n+1),B_{0}(n+i+1))=1\ \  \hbox{for}\ 0\leq
i\leq N,
\end{equation*}
i.e.,
\begin{equation}\label{eq316}
\mathrm{gcd}(A_{0}(n+1),B_{0}(n+i))=1\ \  \hbox{for}\ 1\leq i\leq
N+1,
\end{equation}
and by (\ref{eq38}),
\begin{equation*}
\mathrm{gcd}(A_{0}(n+1),B_{0}(n+k+1))=1\ \  \hbox{for\ all}\
k>N,\end{equation*} i.e.,
\begin{equation}\label{eq317}
\mathrm{gcd}(A_{0}(n+1),B_{0}(n+k))=1\ \  \hbox{for\ all}\ k>N+1.
\end{equation}

\noindent Finally, combining (\ref{eq315}), (\ref{eq316}) and
(\ref{eq317}) into one condition results in
\begin{equation}\label{eq318}
\mathrm{gcd}(A_{0}(n+1),B_{0}(n+h))=1\ \  \hbox{for\ all}\ h\geq 0.
\end{equation}
Hence the right hand side of (\ref{eq314})  is a Gosper
representation for $a(n)/b(n)$.

\begin{ex}\label{ex31} Let $t_{n}=\frac{(n+2)}{n!}$, then
$$r(n)=\frac{t_{n+1}}{t_{n}}=\frac{a(n)}{b(n)},$$ where $a(n)=n+3$,
$b(n)=(n+1)(n+2)$. From Abramov's algorithm, we have $N=1$. By
using Abramov's algorithm with $A_{N+1}(n)=A_{2}(n)=a(n-1)=n+2$
and $B_{N+1}(n)=B_{2}(n)=b(n)=(n+1)(n+2)$, we obtain that
$A_{0}(n)=1$, $B_{0}(n)=n+2$, and the universal denominator
$G(n)=(n+1)(n+2)$. Note that
\begin{equation*}
\mathrm{gcd}(G(n+1),B_{0}(n))=n+2.
\end{equation*}
This means that the Gosper representation $(\ref{eq314})$ in general
is not the $\mathrm{GP}$ representation for $a(n)/b(n)$.
\end{ex}

\noindent As explained in Paule \cite{Pau95}, also the GP
algorithm \cite{Pet92} for finding a GP representation of
$r(n)=a(n)/b(n)$ computes a universal denominator for
(\ref{eq310}), namely as follows.

\noindent {\bf Petkov\v{s}ek's GP Algorithm}: Compute $N$ as in
Abramov's algorithm; i.e., $N=k_{0}$ as in (\ref{eq311}). If
$N=-1$ set $u(n)=1$. If $N\geq0$ define
\begin{equation}\label{eq319}
a_{0}(n)=a(n),\ \ b_{0}(n)=b(n),
\end{equation}
and for $i=1$ up to $i=N+1$ do:
\begin{eqnarray}
\delta_{i}(n)&=& \mathrm{gcd}(a_{i-1}(n),b_{i-1}(n+i)),\label{eq320}\\
a_{i}(n)&=&\frac{a_{i-1}(n)}{\delta_{i}(n)}\ \ \ \mathrm{and}\ \ \
b_{i}(n)=\frac{b_{i-1}(n)}{\delta_{i}(n-i)}.\label{eq321}
\end{eqnarray}
This determines a universal denominator by setting
\begin{equation}\label{eq322}
u(n)=[\delta_{1}(n-1)]^{\underline{1}}[\delta_{2}(n-1)]^{\underline{2}}
\cdots[\delta_{N+1}(n-1)]^{\underline{N+1}}\ .
\end{equation}
Note that this algorithm essentially consists in running the loop
in Abramov's algorithm in the REVERSE direction. (Note that also
its initialization is slightly different, namely starting with the
pair $(a(n),b(n))$ instead of $(a(n-1),b(n))).$

\noindent Finally, as above, from (\ref{eq321}) and by using
(\ref{eq322}) one obtains that
\begin{equation}\label{eq323}
\frac{a(n)}{b(n)}=\frac{u(n+1)}{u(n)}\
\frac{a_{N+1}(n)}{b_{N+1}(n)}.
\end{equation}
Petkov\v{s}ek's algorithm is designed in such a way that
(\ref{eq323}) is not only a Gosper representation, but also a GP
representation for $a(n)/b(n)$; see \cite{Pet92}, \cite{PWZ96},
and also \cite{Pau95}. This means that besides
\begin{equation}\label{eq324}
\mathrm{gcd}(a_{N+1}(n),b_{N+1}(n+h))=1,\ \ \ \hbox{for\ all}\
h\geq0.
\end{equation}
we also have that
\begin{equation}\label{eq325}
\mathrm{gcd}(u(n+1),b_{N+1}(n))=1,
\end{equation}
and
\begin{equation}\label{eq326}
\mathrm{gcd}(u(n),a_{N+1}(n))=1.
\end{equation}
From (\ref{eq314}) and (\ref{eq323}) we get
\begin{equation}\label{eq327}
\frac{u(n+1)}{u(n)}\
\frac{a_{N+1}(n)}{b_{N+1}(n)}=-\frac{G(n+1)}{G(n)}\
\frac{A_{0}(n+1)}{B_{0}(n)}.
\end{equation}
By Lemma 5.3.1 (see \cite{PWZ96}, p.82), we obtain that
$$u(n)\ |\ G(n),$$
which implies that the universal denominator given by GP algorithm
is a factor of the universal denominator given by Abramov's
algorithm.\\

\section{Rational Solutions of Linear Difference Equations}\label{sec4}

In this section we show how to deduce  the explicit formula
(\ref{eq37}) for Abramov's universal denominator by using the
convergence argument.

\begin{thm}\label{thm41}
Let $p_{0}(n)$ and $p_{d}(n)$ be two nonzero polynomials in $n$, and
let $N$ be defined as in $(\ref{eq32})$. Put
\begin{equation}\label{eq41}
G_{k}(n)=\mathrm{gcd}\left(\prod_{j=0}^{k-1}p_{0}(n+j),
\prod_{j=0}^{k-1}p_{d}(n-d-j)\right).
\end{equation}
Then the sequence $G_1(n), G_2(n), \ldots $ converges to
$G_{N+1}(n)$.
\end{thm}

\noindent {\it Proof.} For all $k>N$ we have
\begin{equation}\label{eq42}
\mathrm{gcd}(p_{0}(n), p_{d}(n-d-k)) =\cdots =
\mathrm{gcd}(p_{0}(n+k), p_{d}(n-d-k)) =1.
\end{equation}
Note that
\begin{equation*}
G_{k+1}(n)=\mathrm{gcd}\left(\displaystyle\prod_{j=0}^{k}p_{0}(n+j),
\displaystyle\prod_{j=0}^{k}p_{d}(n-d-j)\right).
\end{equation*}
This implies that
\[ G_k(n) =
G_{k+1}(n), \] for all $k>N$. Moreover, one sees that once
(\ref{eq42}) is satisfied for $k>N$, it is also satisfied for
$k+1$. It follows that
\[ G_k(n)= G_{k+1}(n) = G_{k+2}(n) = \cdots ,\]
for all $k>N$, and this completes the proof. \qed

By using the convergence property of the sequence $G_{k}(n)$, we
can obtain the explicit formula (\ref{eq37}) for Abramov's
universal denominator. Our proof requires only lcm and gcd
computations.

\begin{thm}\label{thm42} Given a linear difference equation
\begin{equation}\label{eq43}
\sum_{m=0}^{d}p_{m}(n)y(n+m)=p(n),
\end{equation}
where $p_{0}(n), p_{1}(n), \ldots, p_{d}(n), p(n)\in\mathbb K[n]$
are given polynomials such that $p_{0}(n)\neq0$, $p_{d}(n)\neq0$.
Let  $f(n)/g(n)$ be the reduced rational representation of $y(n)$,
and $N$ be the dispersion of $p_d(n-d)$ and $p_0(n)$. Then
\begin{equation*}
G(n)=\mathrm{gcd}([p_{0}(n+N)]^{\underline{N+1}},[p_{d}(n-d)]^
{\underline{N+1}}),
\end{equation*}
as given by $(\ref{eq37})$, is a universal denominator of rational
solutions of $(\ref{eq43})$.
\end{thm}

\noindent
{\it Proof.} From (\ref{eq43}) it follows that
\begin{equation}\label{eq44}
\sum_{m=0}^{d}p_{m}(n)\frac{f(n+m)}{g(n+m)}=p(n).
\end{equation}
Letting
\begin{equation}\label{eq45}
l(n)=\mathrm{lcm}(g(n+1),g(n+2),\ldots,g(n+d)),
\end{equation}
and multiplying equation (\ref{eq44}) by $l(n)g(n)$, we obtain
\begin{equation}\label{eq46}
\sum_{m=0}^{d}p_{m}(n)f(n+m)\frac{l(n)}{g(n+m)}\
g(n)=p(n)l(n)g(n).
\end{equation}
From (\ref{eq45}), we have the following divisibility conditions:
\[ g(n+m)\ |\ l(n) \quad \mbox{ for } \; m=1,2,\ldots,d.\]
Thus $l(n)/g(n+m)$ are polynomials for $m=1,2,\ldots,d$. From
(\ref{eq46}) we obtain
\begin{equation}\label{eq47} g(n)\mid
p_{0}(n)\cdot\mathrm{lcm}(g(n+1),g(n+2),\ldots,g(n+d)).
\end{equation}
Similarly, multiplying equation (\ref{eq44}) by $l(n-1)g(n+d)$ and
then substituting $n-d$ for $n$, we obtain that
\begin{equation}\label{eq48} g(n)\mid
p_{d}(n-d)\cdot\mathrm{lcm}(g(n-1),g(n-2),\ldots,g(n-d)).
\end{equation}
Shifting $n$ by 1 in (\ref{eq47}) yields
\begin{equation}\label{eq49}
g(n+1)\mid
p_{0}(n+1)\cdot\mathrm{lcm}(g(n+2),g(n+3),\ldots,g(n+d+1)).
\end{equation}
Substituting (\ref{eq49}) into (\ref{eq47}) we see that $g(n)$
divides
\begin{equation*}
p_{0}(n)\cdot\mathrm{lcm}(p_{0}(n+1)\cdot\mathrm{lcm}(g(n+2),g(n+3),
\ldots,g(n+d+1)),g(n+2),\ldots,g(n+d)).
\end{equation*}
So we can write
\begin{equation*} g(n)\mid
p_{0}(n)p_{0}(n+1)\cdot\mathrm{lcm}(g(n+2),g(n+3),\ldots,g(n+d+1)).
\end{equation*}
By induction we may derive for $k\geq1$,
\begin{equation*}
g(n)\Big|
\prod_{j=0}^{k-1}p_{0}(n+j)\cdot\mathrm{lcm}(g(n+k),g(n+k+1),\ldots,g(n+k+d-1)).
\end{equation*}
It follows that
\begin{equation*}
g(n)\Big|  \prod_{j=0}^{k-1}p_{0}(n+j)\cdot g(n+k)g(n+k+1)\ldots
g(n+k+d-1).
\end{equation*}
Since $\mathbb K$ has characteristic zero, there is a large enough
$k$ such that for any $j\geq k$
$$\mathrm{gcd}(g(n),g(n+j))=1.$$  It follows that
\begin{equation*}
g(n)\Big| \prod_{j=0}^{k-1}p_{0}(n+j),
\end{equation*}
for all large enough $k$. Analogously, from (\ref{eq48}) we get
\begin{equation*}
g(n)\Big| \prod_{j=0}^{k-1}p_{d}(n-d-j),
\end{equation*}
for all large enough $k$. Therefore
\begin{equation*}
g(n)\ |\ G_{k}(n),
\end{equation*}
for all large enough $k$, where $G_{k}(n)$ is defined as in
$(\ref{eq41})$. Setting $k$  to infinity in this equation,  by
Theorem \ref{thm41} we get
\begin{equation*}
g(n)\ |\
G_{N+1}(n)=\mathrm{gcd}\left(\displaystyle\prod_{j=0}^{N}p_{0}(n+j),
\prod_{j=0}^{N}p_{d}(n-d-j)\right)=G(n),
\end{equation*}
as desired. \qed

\vskip 0.5 cm From equation (\ref{eq44}) we get

\begin{equation}\label{eq410}
\sum_{m=0}^{d}p_{m}(n)\cdot f(n+m)\cdot\prod_{\substack{j=0 \\
j\neq m}}^{d}g(n+j)=p(n)\cdot \prod_{j=0}^{d}g(n+j).
\end{equation}
The next step is simply to set
\begin{equation*}
g(n)=G(n),
\end{equation*}
in equation (\ref{eq410}). If equation $\mathrm{(\ref{eq410})}$
can be solved for $f\in\mathbb K[n]$ then $y(n)=f(n)/g(n)$ is a
solution of $\mathrm{(\ref{eq43})}$; otherwise
$\mathrm{(\ref{eq43})}$ has no rational solutions.\\

\begin{alg}\label{alg41}\hskip 0.5 cm \\
\noindent INPUT:  nonzero polynomials
$p_{0}(n), p_{d}(n)$.\\
OUTPUT: a universal denominator $g(n)$ of $(\ref{eq43})$.
\begin{itemize}
\item[$\mathrm {(1)}$] Compute
$N=\mathrm{dis}(p_{d}(n-d),p_{0}(n))=\mathrm{max}\{k\in\mathbb N\
|\ \mathrm{deg\ gcd}(p_{d}(n-d),p_{0}(n+k))\geq1\}.$
\item[$\mathrm {(2)}$] If $N\geq0$ then compute $g(n)=G(n)$, where
$G(n)$ is defined as in $(\ref{eq37})$, otherwise $g(n)=1$.
\end{itemize}
\end{alg}

\begin{ex}\label{ex41}
Find a rational solution of the equation
\begin{eqnarray}\label{eq411}
&&(n+4)(2n+1)(n+2)y(n+3)-(2n+3)(n+3)(n+1)y(n+2)\nonumber\\&&+
n(n+2)(2n-3)y(n+1)-(n-1)(2n-1)(n+1)y(n)=0.
\end{eqnarray}
We have $p_{0}(n)=-(n-1)(2n-1)(n+1)$, $p_{d}(n)=(n+4)(2n+1)(n+2)$,
then $N=2$ and then $g(n)=(n-1)(n+1)n$. By $(\ref{eq410})$, $f(n)$
is a polynomial which satisfies
\begin{eqnarray*}
&&n(2n+1)(n+2)(n+1)f(n+3)-n(2n+3)(n+3)(n+1)f(n+2)+
n(2n-3)\nonumber\\&&\cdot(n+3)(n+2)f(n+1)-(2n-1)(n+3)(n+1)(n+2)f(n)=0.
\end{eqnarray*}
The polynomial $f(n)=Cn(2n-3)$ is a solution of this equation. Thus
\[ y(n)=\frac{f(n)}{g(n)}=C\frac{(2n+1)}{(n^2-1)}\] is a rational
solution of $(\ref{eq411})$.
\end{ex}

\vskip 10pt

\noindent {\bf Acknowledgments.} We would like to thank Qing-Hu Hou
and Doron Zeilberger for helpful discussions and suggestions. This
work was supported by the 973 Project, the PCSIRT Project of the
Ministry of Education, the Ministry of Science and Technology, and
the National Science Foundation of China. Peter Paule was partially
supported by the SFB grant F1305 of the Austrian Science Foundation
FWF.

\end{document}